\documentclass[11pt,a4paper]{article}
\usepackage{epsf, graphicx}
\usepackage{latexsym,amsfonts,amsbsy,amssymb}
\usepackage{amsmath,amsthm}
\usepackage{cite}

\makeatletter

\@addtoreset{equation}{section} \makeatother
\newtheorem{theorem}{Theorem}[section]

\makeatletter \@addtoreset{equation}{section} \makeatother
\renewcommand{\baselinestretch}{1.5}
\begin{document}

\title{\bf A posteriori error estimators suitable for moving finite element methods under anisotropic meshes}
\date{}
 \author{Xiaobo Yin\footnote{Department of
Mathematics, China University of Mining and Technology, Beijing
100083, China (email: yinxb@lsec.cc.ac.cn)} \quad Hehu
Xie\footnote{LSEC, ICMSEC, Academy of Mathematics and Systems
Science, CAS, Beijing 100080, China (email: hhxie@lsec.cc.ac.cn)}}
 \maketitle
 \begin{quote}
\begin{small}
{\bf Abstract.}\,\,In this paper, we give a new type of a posteriori
error estimators suitable for moving finite element methods under
anisotropic meshes for general second-order elliptic problems. The
computation of estimators is simple once corresponding Hessian
matrix is recovered. Wonderful efficiency indices are shown in
numerical experiments.

{\bf Keywords.} a posteriori error estimator; moving finite element
method; anisotropic mesh.

{\bf AMS subject classification.} 65N15, 65N30
\end{small}
\end{quote}
\section{Introduction}
Nowadays adaptive algorithms have been an indispensable tool for
most finite element simulations. They basically consist of the
ingredients ``Solve -- Estimate error -- Refine mesh" which are
repeated until the desired accuracy is achieved. Generally, they can
be classified into three types: $h$-, $r$- and $hp$-version. In this
paper we consider the second ingredient(Estimate error) for
$r$-version(or moving finite element method) under anisotropic
meshes.

Then, what does ``anisotropic mesh" mean? Denote by $h_K$ the
diameter of the finite element $K$, and by $\varrho_K$ the supremum
of the diameters of all balls contained in $K$. It is assumed in the
classical finite element theory that
\begin{eqnarray}\label{regular}
h_K\lesssim \varrho_K.
\end{eqnarray}
(The notation $\lesssim$ means smaller than up to a constant.)
Elements which satisfy (\ref{regular}) are called isotropic
elements.

Many physical problems exhibit a common anisotropic feature that
their solutions change more significantly in one direction than the
others. Examples include those having boundary layers, shock waves,
interfaces, and edge singularities, etc.. In such cases it is
advantageous to reflect this anisotropy in the discretization by
using meshes with anisotropic elements (sometimes also called
elongated elements). These elements have a small mesh size in the
direction of the rapid variation of the solution and a larger mesh
size in the perpendicular direction. That is to say, anisotropic
elements do not satisfy condition (\ref{regular}). Conversely they
are characterized by
\begin{eqnarray}\label{notregular}
\frac{h_K}{\varrho_K}\rightarrow \infty
\end{eqnarray}
where the limit can be considered as $h\rightarrow 0$ (near edges)
or $\epsilon \rightarrow 0$ (in layers) where $\epsilon$ is some
(small perturbation) parameter of the problem. Indeed anisotropic
meshes have been used successfully in many areas, for example in
singular perturbation and flow problems
\cite{Ait-Ali-Yahia,ApelLube,Becker,HabForDomValBou,PerVahMorZ,ZWu}
and in adaptive procedures
\cite{BuscagliaDari,CasHecMohPir,PerVahMorZ,Siebert}. For problems
with very different length scales in different spatial directions,
long and thin triangles turn out to be better choices than shape
regular ones if they are properly used. This motivated an intensive
study on the error analysis for anisotropic meshes in the finite
element method. For instance, Apel \cite{Apel} described an error
estimate in terms of the length scales $h_1$ and $h_2$ along the $x$
and $y$ direction, respectively. Berzins \cite{Berzins} developed a
mesh quality indicator measuring the correlation between the
anisotropic features of the mesh and those of the solutions. Kunert
\cite{Kunert2000} introduced the concept of ``matching function"
which measures the correspondence between an anisotropic mesh and a
given function. Using this concept he gave three types of a
posteriori error estimator for anisotropic meshes under the
assumption that the anisotropic mesh $T_h$ is `adapted' to the
anisotropic solution. Formaggia and Perotto \cite{ForPer} used the
spectral properties of the affine mapping from a reference triangle
to obtain a full information about the orientation, dimension and
aspect ratio of a given element. After that they proposed a
posteriori estimators for elliptic problems under anisotropic
meshes. Picasso\cite{Picasso} combined the method in \cite{ForPer}
and a Zienkiewicz-Zhu error estimator to approach the error
gradient. Cao \cite{cao} revealed the precice relation between the
error of linear interpolation on a general triangle and the
geometric characters of the triangle. This list is certainly
incomplete, but from the papers we can find the interpolation error
depends on the solution and the size and shape of the elements in
the mesh.

In the mesh generation community, the error estimate is often
studied for the model problem of interpolating quadratic functions.
This model is a reasonable simplification of the cases involving
general functions, since quadratic functions are the leading terms
in the local expansion of the linear interpolation errors. For
instance, Nadler \cite{Nadler} derived an exact expression for the
$L^2$-norm of the linear interpolation error in terms of the three
sides ${\bf \ell}_1$, ${\bf \ell}_2$, and ${\bf \ell}_3$ of the
triangle $K$,
\begin{eqnarray}\label{Nadlerformula}
\|u-u_I\|^2_{L^2(K)}=\frac{|K|}{180}{\Big
[}(d_1+d_2+d_3)^2+d_1d_2+d_2d_3+d_1d_3{\Big ]},
\end{eqnarray}
where $|K|$ is the area of the triangle, $d_i = {\bf \ell}_i\cdot
H{\bf \ell}_i$ with $H$ being the Hessian matrix of $u$. Bank and
Smith \cite{BankSmith} gave a formula for the $H^1$-seminorm of the
linear interpolation error
\begin{eqnarray}\label{Bankformula}
\|\nabla(u-u_I)\|^2_{L^2(K)}=\frac{1}{4}{\bf d}\cdot B{\bf d},
\end{eqnarray}
where ${\bf d}=[d_1,d_2,d_3]^T$,
\begin{eqnarray*}
B=\frac{1}{48|K|}\left(
\begin{array}{ccc}
 |{\bf \ell}_1|^2+|{\bf \ell}_2|^2+|{\bf \ell}_3|^2 & 2{\bf \ell}_1\cdot{\bf \ell}_2 & 2{\bf \ell}_1\cdot{\bf \ell}_3 \\
 2{\bf \ell}_1\cdot{\bf \ell}_2 & |{\bf \ell}_1|^2+|{\bf \ell}_2|^2+|{\bf \ell}_3|^2  & 2{\bf \ell}_2\cdot{\bf \ell}_3\\
 2{\bf \ell}_1\cdot{\bf \ell}_3 & 2{\bf \ell}_2\cdot{\bf \ell}_3 &|{\bf \ell}_1|^2+|{\bf \ell}_2|^2+|{\bf \ell}_3|^2
\end{array}
\right)
\end{eqnarray*}
In this paper we'll develop the formula for $H^1$-seminorm of the
linear interpolation error
\begin{eqnarray}
\|\nabla(u-u_I)\|^2_{L^2(K)}\approx\sum_{K\in
T_h}\frac{1}{48|K|}\sum_{i=1}^{3}c_i^2|{\bf \ell}_i|^2,
\end{eqnarray}
where $c_i={\bf \ell}_{i+1}\cdot H{\bf \ell}_{i+2}$, and that of
discretization error
\begin{eqnarray}
\|\nabla(u-u_h)\|^2_{L^2(\Omega)}\approx-\frac{1}{24}\sum_{K\in
T_h}\sum_{i=1}^{3}{\Big(}f_K+|{\bf \ell}_i|[\partial_{n}u_h]_{{\bf
\ell}_i}{\Big)}d_i.
\end{eqnarray}

The quality of an a posteriori error estimator is often measured by
its efficiency index, i.e., the ratio of the true error and the
estimated error(in some norm). An error estimator is called
efficient if its efficiency index together with its inverse remain
bounded for all mesh-sizes. It is called asymptotically exact if its
efficiency index tends to one when the mesh-size converges to zero.
From numerical results we see our estimators are often
asymptotically exact although we couldn't prove it rigorously.

The paper is organized as follows. In section 2 we give some
preliminary results, especially the error expansions for $u-u_I$ and
$\nabla(u-u_I)$. In section 3 these error expansions are used to
derive a posteriori error estimators for the interpolation error and
discretization error, respectively.  Section 4 contains ``efficient
index" tables and pictures from numerical experiments for some
second-order elliptic problems which yield anisotropic solutions.
The results show remarkable agreement with the theoretical
predictions. Finally, in section 5 we state our conclusions and
direction for further research.

\section{Preliminaries}
Consider the following model problem. Find $u$: $\Omega \subset
\mathcal{R}^2\rightarrow \mathcal{R}$ such that
\begin{equation}\label{2.01}
     \left \{
     \begin{array}{lll}
     Lu&=&-\sum\limits_{i,j=1}^{2}\frac{\partial }{\partial x_i}{\Big (}a_{ij}\frac{\partial u}{\partial x_j}{\Big
     )}
     +bu=f\quad{\rm in} \,\, \Omega,\\
     u&=&0\quad{\rm on}\,\,  \partial\Omega,
     \end{array}
     \right .
\end{equation}
where $b=b({\bf x})\geq 0$ a.e. in $\Omega$ and $a_{ij}=a_{ij}({\bf
x})$ are given functions. The domain $\Omega$ is an open, bounded
subset of $\mathcal{R}^2$ and the operator $L$ is elliptic and
self-adjoint. The corresponding variational formulation seeks $u\in
H_0^1(\Omega)$ such that
\begin{eqnarray}\label{2.02}
a(u,v)=(f,v) \quad \forall v\in V\equiv H_0^1(\Omega),
\end{eqnarray}
where
\begin{eqnarray*}
a(u,v)\equiv\int_{\Omega}{\Big
(}\sum\limits_{i,j=1}^{2}a_{ij}\frac{\partial u}{\partial
x_j}\frac{\partial v}{\partial x_i}+buv{\Big )}d{\bf x}
\end{eqnarray*}
and
\begin{eqnarray*}
(f,v)\equiv\int_{\Omega}fvd{\bf x}.
\end{eqnarray*}
We shall use the standard notations in \cite{Ciarlet} for the
Sobolev spaces $H^s(\Omega)$ and their associated inner products
$(\cdot,\cdot)_s$, norms $||\cdot||_s$, and seminorms $|\cdot|_s$
for $s \geq 0$.

By $\mathcal{F}=\{\mathcal{T}_h\}$ we denote a family of
triangulations $\mathcal{T}_h$ of $\Omega$. Let $V_h$ be the space
of continuous, piecewise linear functions over $\mathcal{T}_h$, and
$V_{0,h}\equiv V_h\cap H_0^1(\Omega)$. The finite element
approximation problem of (\ref{2.02}) seeks $u_h\in V_{0,h}$ such
that
\begin{eqnarray}\label{2.03}
a(u_h,v_h)=(f,v_h) \quad \forall v_h\in V_{0,h}.
\end{eqnarray}
The three vertices of an arbitrary triangle $K\in \mathcal{T}_h$ are
denoted by ${\bf a_1}=(x_1,y_1)^T$, ${\bf a_2}=(x_2,y_2)^T$ and
${\bf a_3}=(x_3,y_3)^T$. Additionally we define the edge vectors
${\bf \ell_1}={\bf a_3}-{\bf a_2}$, ${\bf \ell_2}={\bf a_1}-{\bf
a_3}$ and ${\bf \ell_3}={\bf a_2}-{\bf a_1}$(Figure 1).
\begin{figure}[ht!]
  \centering
  \includegraphics[width=4cm]{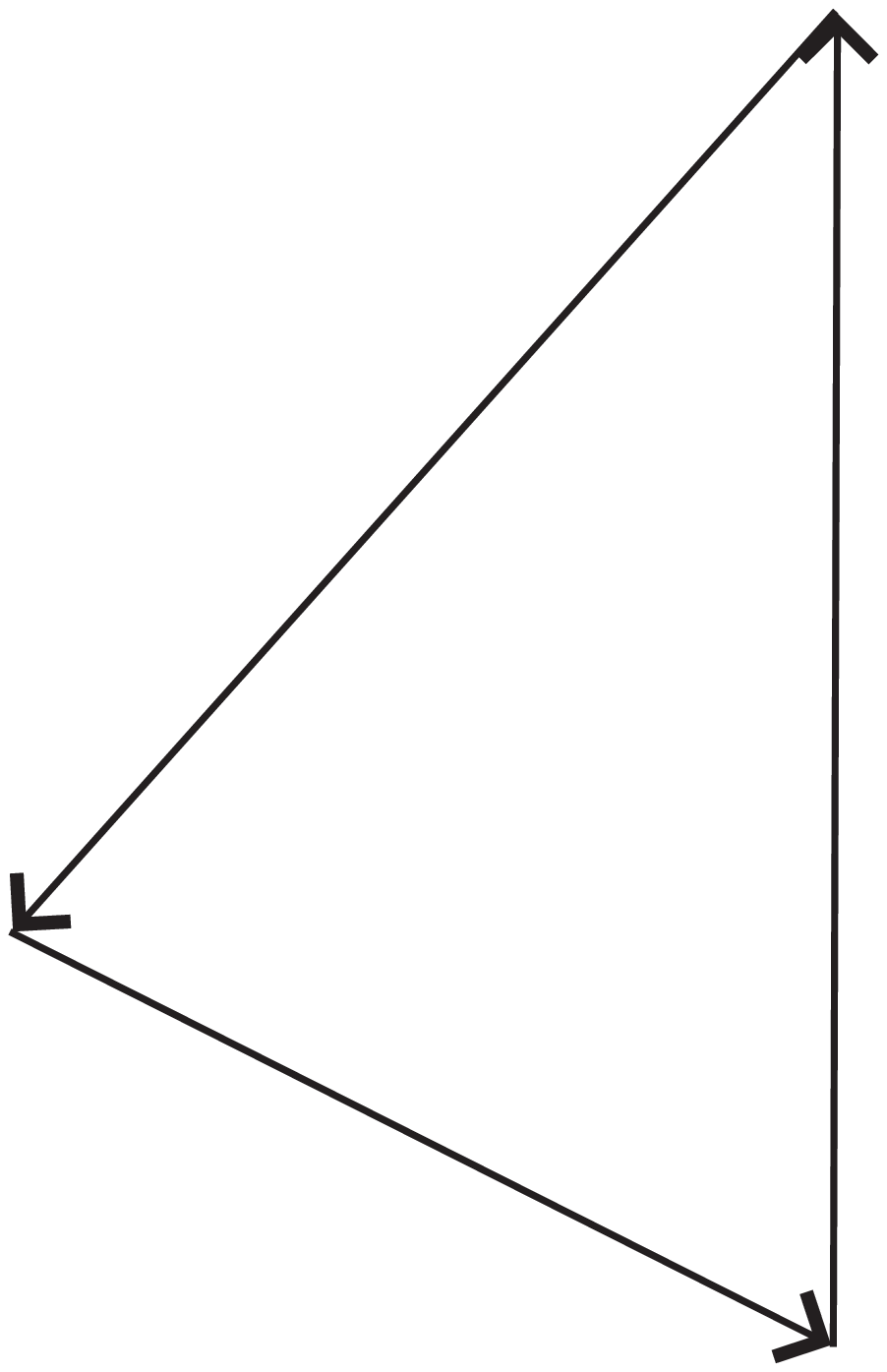}
 \put(0,160){\mbox{${\bf a}_1$}}
\put(-125,55){\mbox{${\bf a}_2$}} \put(-5,0){\mbox{${\bf a}_3$}}
\put(-50,80){\mbox{$K$}} \put(-70,20){\mbox{${\bf \ell}_1$}}
 \put(-5,80){\mbox{${\bf
\ell}_2$}} \put(-75,110){\mbox{${\bf \ell}_3$}}
\put(-160,-20){\mbox{Figure 1: notations in a single element $K$.}}
\end{figure}
Denote by $u_I$ the linear interpolation of $u$ at the three
vertices of $K$. Let $\{\lambda_i({\bf x})\}_{i=1}^3$ be the
barycentric coordinates of $K$. From \cite{ChenSunXu,Shewchuk} we
know for a quadratic function $u$ over $K$ the following formulas
hold:
\begin{eqnarray}\label{2.1_L}
u({\bf x})-u_I({\bf
x})=-\frac{1}{2}\sum\limits_{i=1}^{3}\lambda_i({\bf x}){\Big [}({\bf
x}-{\bf a_i})\cdot H({\bf x}-{\bf a_i}){\Big ]}\,\,\forall {\bf
x}\in K,
\end{eqnarray}
\begin{eqnarray}\label{2.1}
\nabla{\Big (}u({\bf x})-u_I({\bf x}){\Big
)}=-\frac{1}{2}\sum\limits_{i=1}^{3}\nabla\lambda_i({\bf x}){\Big
[}({\bf x}-{\bf a_i})\cdot H({\bf x}-{\bf a_i}){\Big ]}\,\, \forall
{\bf x}\in K,
\end{eqnarray}
where $H$ is the Hessian matrix of $u$.

\section{A posteriori error estimates}
Many authors have discussed the interpolation error to derive their
adaptive algorithm(\cite{Apel,BankSmith,Berzins,ForPer,Huang}).
However, the interpolation error is different from the
discretization error in most cases. In this section we first discuss
the former and then the latter. Finally we will analyze their
relationship using the concept ``superapproximation''.
\subsection{An a posteriori error estimator for the interpolation error}
\begin{theorem}\label{Theorem3.1}
Let $u$ be a quadratic function and $u_I$ is the Lagrangian linear
finite element interpolation of $u$. Denote by $H$ the Hessian
matrix of $u$. The following relationship holds:
\begin{eqnarray}\label{2.2}
\|\nabla(u-u_I)\|^2_{L^2(K)}=\frac{1}{48|K|}\sum_{i=1}^{3}c_i^2|{\bf
\ell}_i|^2,
\end{eqnarray}
where $c_i={\bf \ell}_{i+1}\cdot H{\bf \ell}_{i+2}$. Here we
prescribe $i+3=i,i-3=i$.
\end{theorem}

\begin{proof}
From (\ref{2.1}), we have
\begin{eqnarray}
&&\|\nabla(u-u_I)\|^2_{L^2(K)}=\frac{1}{4}\int_K{\Big|}\sum\limits_{i=1}^{3}
\nabla\lambda_i({\bf x}){\Big [}({\bf x}-{\bf a}_i)\cdot H({\bf
x}-{\bf a}_i){\Big ]}{\Big |}^2d{\bf x}.
\end{eqnarray}
Due to the properties of the barycentric coordinates it is known
that
\begin{eqnarray}
{\Big|}\sum\limits_{i=1}^{3} \nabla\lambda_i({\bf x}){\Big [}({\bf
x}-{\bf a}_i)\cdot H({\bf x}-{\bf a}_i){\Big ]}{\Big |}^2\in P_2(K).
\end{eqnarray}
We use the second-order quadrature scheme which is exact for
polynomial of degree less or equal to 2, i.e.
\begin{eqnarray}\label{2.3}
\int_K\varphi({\bf x})d{\bf
x}=\frac{1}{3}{\Big|}K{\Big|}{\Big[}\varphi({\bf
a}_{12})+\varphi({\bf a}_{23})+\varphi({\bf a}_{31}){\Big
]},\quad\forall \varphi\in P_2(K),
\end{eqnarray}
where ${\bf a}_{ij}$ is the midpoint of the segment $\overline{{\bf
a}_i{\bf a}_{j}}$. Notice that
\begin{eqnarray*}
\nabla\lambda_i({\bf
x})={\Big(}\frac{y_{i+1}-y_{i+2}}{2|K|},-\frac{x_{i+1}-x_{i+2}}{2|K|}{\Big)}^T,
\end{eqnarray*} after a simple calculation we get
(\ref{2.2}).
\end{proof}

Here we set
\begin{eqnarray}\label{2.4}
\eta_I=\sqrt{\sum_{K\in T_h}\frac{1}{48|K|}\sum_{i=1}^{3}c_i^2|{\bf
\ell}_i|^2}
\end{eqnarray}
as the a posteriori estimator for $\|\nabla(u-u_I)\|_{L^2(\Omega)}$. \\
{\bf Remark 1} Using the same technique we can also get the
corresponding estimator for $\|u-u_I\|_{0,\Omega}$ denoted by
$\eta_{I0}$.\\
{\bf Remark 2} The estimators $\eta_I$ and $\eta_{I0}$ have been
given in different forms, for example, in
\cite{BankSmith,cao,Nadler}, e.t.c..

\subsection{An a posteriori error estimator for the discretization error}
In this subsection an a posteriori error estimator for the
discretization error of problem (\ref{2.02}) will be given.
\begin{theorem}\label{Theorem3.2}
Assume $L$ be an elliptic and adjoint operator, $c({\bf x})$ and
$a_{ij}({\bf x})$ be zero and constant functions, respectively. We
have the following estimate:
\begin{eqnarray}\label{3.2}
\int_{\Omega}{\Big|}\nabla(u-u_h){\Big|}^2d{\bf
x}\approx-\frac{1}{24}\sum_{K\in T_h}\sum_{i=1}^{3}{\Big(}f_K+|{\bf
\ell}_i|[\partial_{n}u_h]_{{\bf \ell}_i}{\Big)}d_i,
\end{eqnarray}
where $f_K=\int_Kf({\bf x})d{\bf x}$ and $[\partial_{n}u_h]_{{\bf
\ell}_i}$ is the jump of the conormal derivative of $u_h$
$$\partial_{n}u_h=\sum_{i,j=1}^{2}a_{ij}\frac{\partial u_h}{\partial x_j}n_i$$
across the edge ${\bf \ell}_i$, with ${\bf n}=(n_1,n_2)^T$ the unit
outward normal vector.
\end{theorem}

\begin{proof}
Using the Galerkin orthogonality, we have
\begin{eqnarray*}
&&\int_{\Omega}{\Big|}\nabla(u-u_h){\Big|}^2d{\bf
x}=a(u-u_h,u-u_h)=a(u-u_h,u-u_I)\nonumber\\
&=&\sum_{K\in T_h}\int_K{\Big[}\sum_{i,j=1}^{2}a_{ij}\frac{\partial
(u-u_h)}{\partial
x_j}\frac{\partial}{\partial x_i}(u-u_I){\Big]}d{\bf x}\nonumber\\
&=&\sum_{K\in
T_h}{\Big\{}-\int_K{\Big[}\sum_{i,j=1}^{2}\frac{\partial}{\partial
x_i}{\Big(}a_{ij}\frac{\partial (u-u_h)}{\partial
x_j}{\Big)}(u-u_I){\Big]}d{\bf x}+\int_{\partial
K}\partial_{n}(u-u_h)(u-u_I)ds{\Big\}}\nonumber\\
&=&\sum_{K\in T_h}{\Big\{}\int_K\frac{f_K}{|K|}(u-u_I)d{\bf
x}+\frac{1}{2}\int_{\partial
K}[\partial_{n}u_h](u-u_I)ds{\Big\}}+\sum_{K\in T_h}\int_K{\Big(}f-\frac{f_K}{|K|}{\Big)}(u-u_I)d{\bf x}\nonumber\\
&\approx&\sum_{K\in T_h}{\Big\{}\int_K\frac{f_K}{|K|}(u-u_I)d{\bf
x}+\frac{1}{2}\int_{\partial
K}[\partial_{n}u_h](u-u_I)ds{\Big\}}\nonumber\\
&\approx&-\frac{1}{24}\sum_{K\in
T_h}{\Big(}f_K\sum_{i=1}^{3}d_i+\sum_{i=1}^{3}|{\bf
\ell}_i|[\partial_{n}u_h]_{{\bf
\ell}_i}d_i{\Big)}=-\frac{1}{24}\sum_{K\in
T_h}\sum_{i=1}^{3}{\Big(}f_K+|{\bf \ell}_i|[\partial_{n}u_h]_{{\bf
\ell}_i}{\Big)}d_i,
\end{eqnarray*}
where we use the error expansion (\ref{2.1_L}) and the second-order
quadrature scheme on $K$ and $\partial K$, respectively.
\end{proof}

\subsection{Discussion of the estimators}
From the Theorem \ref{Theorem3.2} we get easily an a posteriori
error estimator for the discretization error:
$$\eta=\sqrt{-\frac{1}{24}\sum_{K\in T_h}\sum_{i=1}^{3}{\Big(}f_K+|{\bf
\ell}_i|[\partial_{n}u_h]_{{\bf \ell}_i}{\Big)}d_i}.$$ Obviously
this estimator can be computed easily provided that $H$ is properly
given. A number of numerical recovery approaches have been proposed
in the literature for second-order
derivatives\cite{Agou,Lip,Zhangxd,Zhangzm,ZZ1987}. Comparisons of
these techniques have also been made in \cite{Buscaglia,Vallet}.
Particularly the authors\cite{Vallet} compared four methods for
reconstructing the second-order derivatives of a piecewise linear
function: DLF(Double linear fitting), SLF(Simple linear fitting),
QF(Quadratic fitting) and DL2P(Double $L^2$-projection). In this
paper we will recover $H$ using the quadratic fitting method
elaborated by Zhang \cite{Zhangxd}.

To end this section, it is advantageous to discuss the relationship
between the interpolation error $\|\nabla(u-u_I)\|_{0,\Omega}$ and
the discretization error $\|\nabla(u-u_h)\|_{0,\Omega}$.

Denote by $N$ the number of elements in $\mathcal{T}_h$. Assume
$\|\nabla(u-u_I)\|_{0,\Omega}\approx CN^{1/2}$ and
$\|\nabla(u-u_h)\|_{0,\Omega}\approx CN^{1/2}$. Then, by simple
calculus, we have
\begin{eqnarray}\label{3.7}
\|\nabla(u-u_I)\|_{0,\Omega}^2-\|\nabla(u-u_h)\|_{0,\Omega}^2=\int_\Omega
\nabla(u_I+u_h-2u)\cdot(\nabla u_I-\nabla u_h)d{\bf x}.
\end{eqnarray}
From (\ref{3.7}) we conclude that if
\begin{eqnarray}\label{3.75}
\|\nabla(u_I-u_h)\|_{0,\Omega}\leq CN^{-\frac{1}{2}-\gamma},
\end{eqnarray}
where $\gamma>0$(this phenomena is called
superapproximation\cite{BanXu1,LinLin}), then
\begin{eqnarray}
\|\nabla(u-u_I)\|_{0,\Omega}^2-\|\nabla(u-u_h)\|_{0,\Omega}^2=O(N^{-1-\gamma}).
\end{eqnarray}
Assume $\eta_I$ be an asymptotically exact estimator of
$\|\nabla(u-u_I)\|_{0,\Omega}$, that is
\begin{eqnarray}\nonumber
\lim\limits_{N\rightarrow
\infty}\frac{\eta_I}{\|\nabla(u-u_I)\|_{0,\Omega}}=1.
\end{eqnarray}
Then $\eta_I$ can also be used as the estimator of
$\|\nabla(u-u_h)\|_{0,\Omega}$ because
\begin{eqnarray}
\lim\limits_{N\rightarrow
\infty}\frac{\eta_I^2}{\|\nabla(u-u_h)\|_{0,\Omega}^2}=\lim\limits_{N\rightarrow
\infty}\frac{\eta_I^2}{\|\nabla(u-u_I)\|_{0,\Omega}^2+O(N^{-1-\gamma})}=1.
\end{eqnarray}
However, the superapproximation can be proved only in some
structured meshes such as uniform and uniform Chevron triangular
meshes in \cite{LinLin}, and $O(h^{2\sigma})$ irregular triangular
meshes in \cite{BanXu1}, under the assumption that $u$ is very
smooth. When the solution doesn't have superapproximation property
we couldn't replace the discretization error by the interpolation
error. Fortunately, from numerical experiments in section 4 we guess
that the superapproximation always holds during the adaptive
procedure.

\subsection{Problem for general coefficients}
For the discussion above we assume that $b({\bf x})$ and
$a_{ij}({\bf x})$ are  zero and constant functions, respectively. In
fact, we can get the corresponding results for the general smooth
functions $b({\bf x})$ and $a_{ij}({\bf x})$, if we use the higher
order quadrature scheme and notice that
\begin{eqnarray*}
&&\sum_{K\in T_h}\int_K{\Big[}\sum_{i,j=1}^{2}a_{ij}({\bf
x})\frac{\partial (u-u_h)}{\partial x_j}\frac{\partial}{\partial
x_i}(u-u_I){\Big]}d{\bf
x}\nonumber\\
&&=\sum_{K\in
T_h}{\Big\{}-\int_K{\Big[}\sum_{i,j=1}^{2}\frac{\partial}{\partial
x_i}{\Big(}a_{ij}({\bf x})\frac{\partial (u-u_h)}{\partial
x_j}{\Big)}(u-u_I){\Big]}d{\bf
x}+\int_{\partial K}\partial_{n_{\bf x}}(u-u_h)(u-u_I)ds{\Big\}}\nonumber\\
&&=\sum_{K\in
T_h}{\Big\{}-\int_K{\Big[}\sum_{i,j=1}^{2}\frac{\partial}{\partial
x_i}{\Big(}a_{ij}^{I}({\bf x})\frac{\partial (u-u_h)}{\partial
x_j}{\Big)}(u-u_I){\Big]}d{\bf
x}+\int_{\partial K}\partial_{n}(u-u_h)(u-u_I)ds{\Big\}}\nonumber\\
&&+\sum_{K\in
T_h}{\Big\{}\int_K{\Big[}\sum_{i,j=1}^{2}\frac{\partial}{\partial
x_i}{\Big(}{\Big(}a_{ij}^{I}({\bf x})-a_{ij}({\bf
x}){\Big)}\frac{\partial
(u-u_h)}{\partial x_j}{\Big)}(u-u_I){\Big]}d{\bf x}\nonumber\\
&&+\int_{\partial K}(\partial_{n_{\bf
x}}-\partial_{n})(u-u_h)(u-u_I)ds{\Big\}}\nonumber\\
&&\approx\sum_{K\in
T_h}{\Big\{}-\int_K{\Big[}\sum_{i,j=1}^{2}\frac{\partial}{\partial
x_i}{\Big(}a_{ij}^{I}({\bf x})\frac{\partial (u-u_h)}{\partial
x_j}{\Big)}(u-u_I){\Big]}d{\bf x}+\int_{\partial
K}\partial_{n}(u-u_h)(u-u_I)ds{\Big\}},
\end{eqnarray*}
where $a_{ij}^{I}({\bf x})$ is the Lagrangian linear finite element
interpolant of $a_{ij}({\bf x})$, and
$$\overline{a}_{ij}({\bf
x})=\frac{1}{|K|}\int_Ka_{ij}({\bf x})d{\bf
x},\quad\partial_{n}u_h=\sum\limits_{i,j=1}^{2}\overline{a}_{ij}({\bf
x})\frac{\partial u_h}{\partial x_j}n_i,\quad\partial_{n_{\bf
x}}u_h=\sum\limits_{i,j=1}^{2}a_{ij}({\bf x})\frac{\partial
u_h}{\partial x_j}n_i.$$

\section{Numerical experiments}
First we give some definitions,
$$E=\frac{\eta^2}{\|\nabla(u-u_h)\|_{0,\Omega}^2}, \quad EI=\frac{\eta_I^2}{\|\nabla(u-u_h)\|_{0,\Omega}^2},\quad
\textrm{Where}\,\, H\,\,\textrm{is}\,\,\textrm{exact}\,\,
\textrm{Hessian}\,\, \textrm{matrix},$$
$$E_r=\frac{\eta^2}{\|\nabla(u-u_h)\|_{0,\Omega}^2}, \quad EI_r=\frac{\eta_I^2}{\|\nabla(u-u_h)\|_{0,\Omega}^2},\quad
\textrm{Where}\,\,H_r\,\,\textrm{is}\,\, \textrm{recovered}\,\,
\textrm{by} \,\, \textrm{Zhang}\mbox{\cite{Zhangxd}}.$$ We induce
the exact Hessian for comparison in examples 4.1-4.4 and 4.6, while
in example 4.5 we just use the recovered $H_r$ where the true
solution doesn't belong to $H^2(\Omega)$.

Because we use the Hessian recovery technique in \cite{Zhangxd}, it
is advantageous to show how this technique works. From this point we
will verify if there exists a positive number $\delta$ such that
$||H-H_r||_{0,\Omega}<CN^{-\delta}$ in our numerical experiment.\\
{\bf Example 4.1} This example is to solve the boundary value
problem of Poisson's equation
\begin{eqnarray}
-\triangle u&=&f,\quad {\bf x}\in \Omega\equiv(0,1)\times(0,1),
\end{eqnarray}
where the Dirichlet boundary condition and the right-hand side term
are chosen such that the exact solution is given by
\begin{eqnarray}
u({\bf x})={\Big(}1+e^{\frac{x_1+x_2-0.85}{2\epsilon}}{\Big)}^{-1}
\end{eqnarray}
with $\epsilon$ being taken to be 0.005(taken from \cite{Huang}).
Here we use the Delauney mesh generator to get the nearly uniform
mesh, where $nu$ is the number of initial points on the boundary.
See Table 1 and Figure 2 for more details.
\begin{figure}[ht]
  \includegraphics[width=8.5cm]{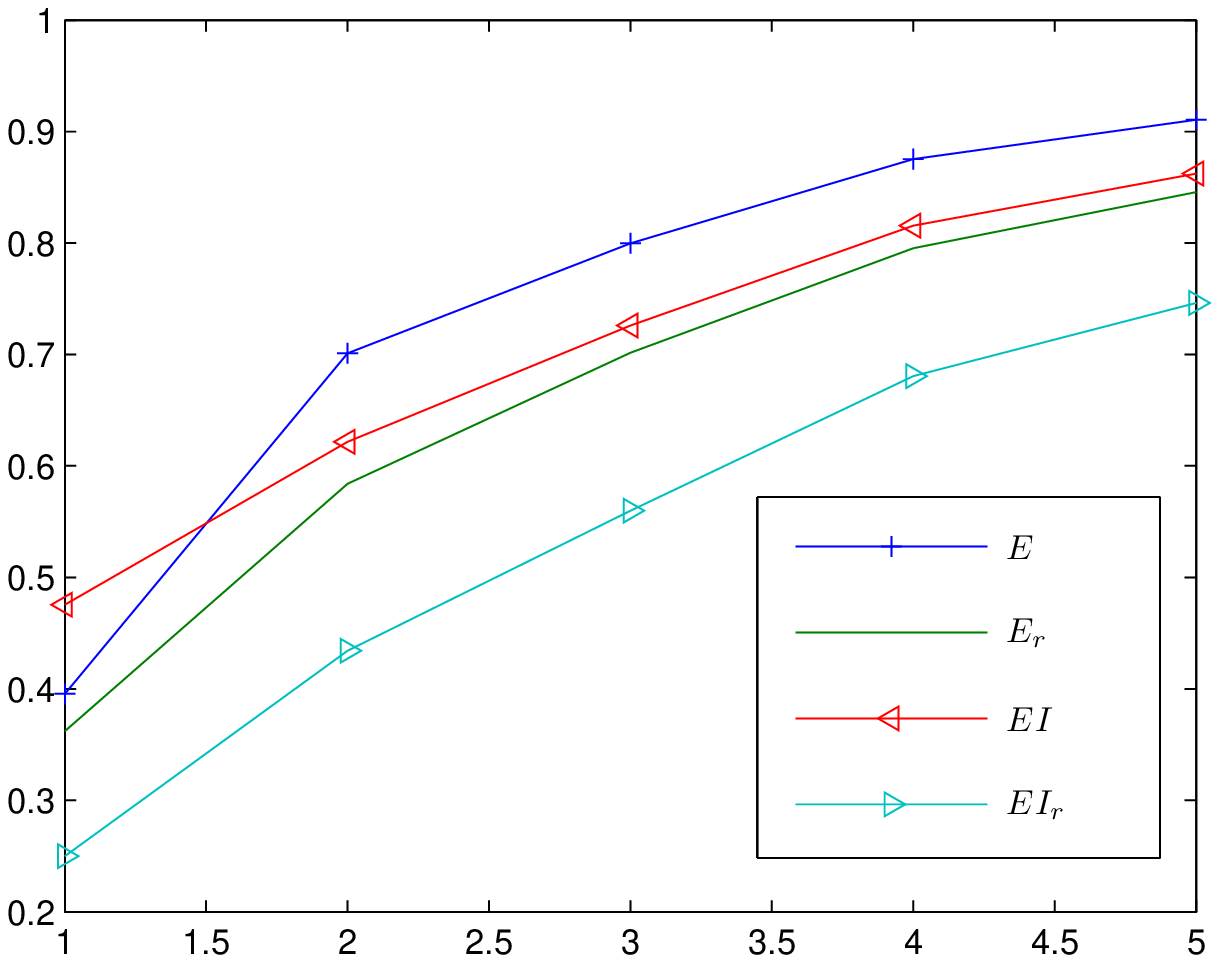}
  \put(-220,0){(a)}
  \put(0,16){\resizebox{6.5cm}{!}{\epsffile{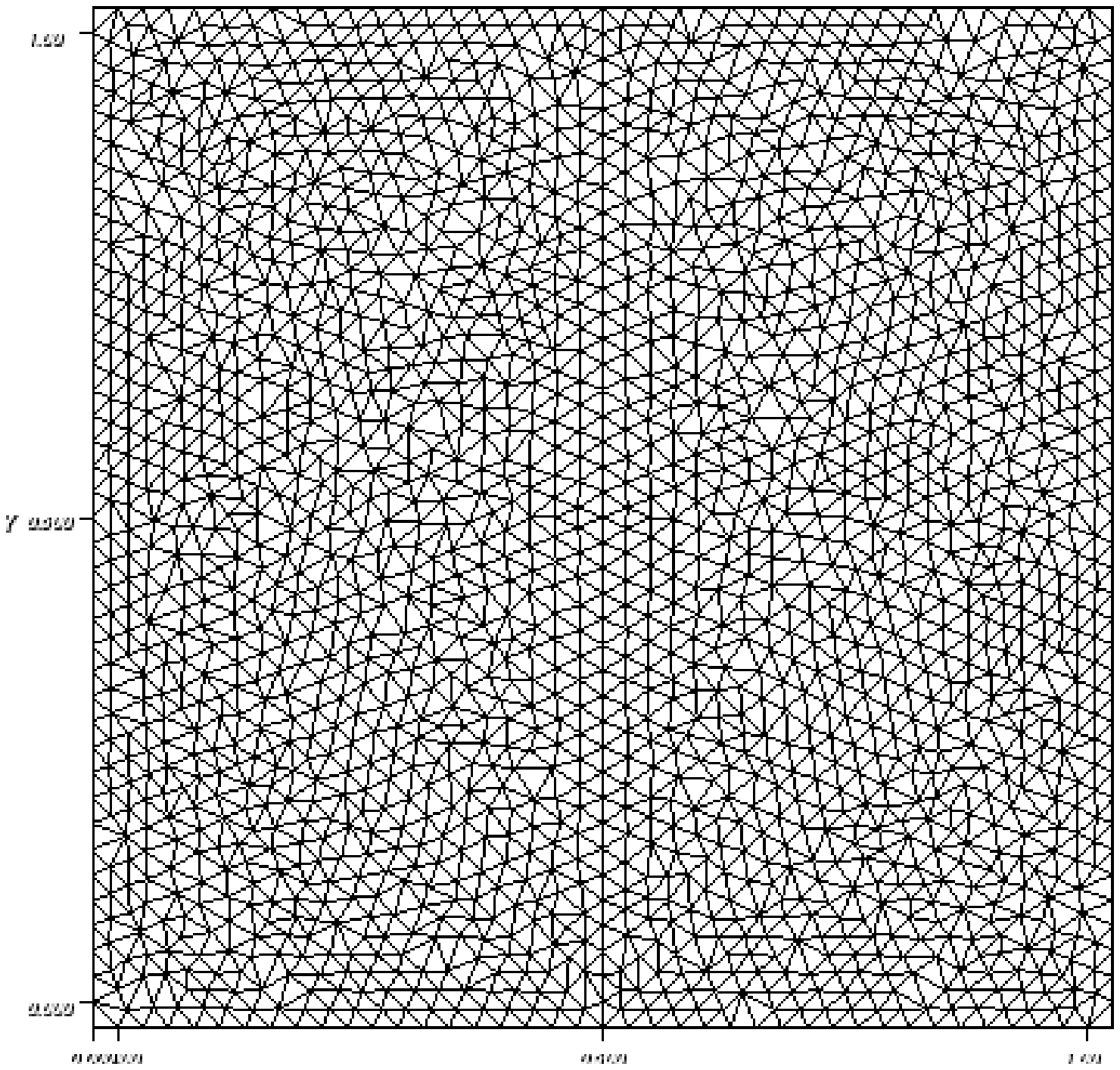}}}
  \put(10,0){(b)}
  \put(-230,-20){\mbox{Figure 2: (a) Four
estimators and (b) the initial mesh of the example 4.1.}}
\end{figure}

\begin{center}
{\renewcommand{\baselinestretch}{1.0}\small
\begin{tabular}{cccccccc}
\multicolumn{8}{c}{\hbox{{\bf Table 1: Four estimators and $\delta$ in example~4.1}~~  }}\\
 \hline
$nu$ & $N$ &$E$ & $E_r$ &$EI$ &$EI_r$ & $||H-H_r||$& $\delta$\\
\hline
40  & 3744  &0.395664 & 0.361983 & 0.475616 & 0.250096 &172.773 & -\\
80  & 8664  &0.701157 & 0.584010  & 0.621573 & 0.434262 &92.8695 & 1.48\\
160 & 15154 &0.799865 & 0.701711 & 0.726035 & 0.560006 &62.6707 & 1.41\\
320 & 23674 &0.875311 & 0.795383 & 0.815560 & 0.680712  &40.1777 & 1.99\\
640 & 34108 &0.910547 & 0.845726 & 0.862312 & 0.746308  &29.0403 & 1.78\\
\hline
\end{tabular}}
\end{center}
From example 4.2 to 4.6 meshes are generated using a c++ code
BAMG(Bidimensional Anisotropic Mesh Generator) developed by Hecht
\cite{Hecht} ($r-$version adaptive procedure). \\
 {\bf Example 4.2} The same problem as in example 4.1.
In fact the solution exhibits a sharp layer on line
$x_1+x_2-0.85=0$. The result is shown in the following table, where
$p$ stands for the step of the adaptive procedure. Results are
listed in Table 2 and Figure 3.

\begin{center}
{\renewcommand{\baselinestretch}{1.0}\small
\begin{tabular}{cccccccc}
\multicolumn{8}{c}{\hbox{{\bf Table 2: Four estimators and $\delta$ in example~4.2}~~  }}\\
 \hline
$p$ & $N$ &$E$ & $E_r$ &$EI$ &$EI_r$ & $||H-H_r||$& $\delta$\\
\hline
1 & 94  &-0.015593 & 0.136566 & 0.812635 & 0.127626 & 265.946  &  -\\
2 & 113 & 0.035375 & 0.176612 & 0.849958& 0.179425 & 233.315 & 1.42\\
3 & 189 &0.668308  & 0.540584 & 0.601845 & 0.414778 & 153.646 & 1.62\\
4 & 272 &0.950449 & 0.846033  & 0.929412 & 0.738008 & 62.0580 & 4.98\\
5 & 278 &0.994030 & 0.916821  & 1.018023 & 0.850781 & 34.0592 & 55.0\\
\hline
\end{tabular}}
\end{center}

\begin{figure}[ht]
  \includegraphics[width=8.5cm]{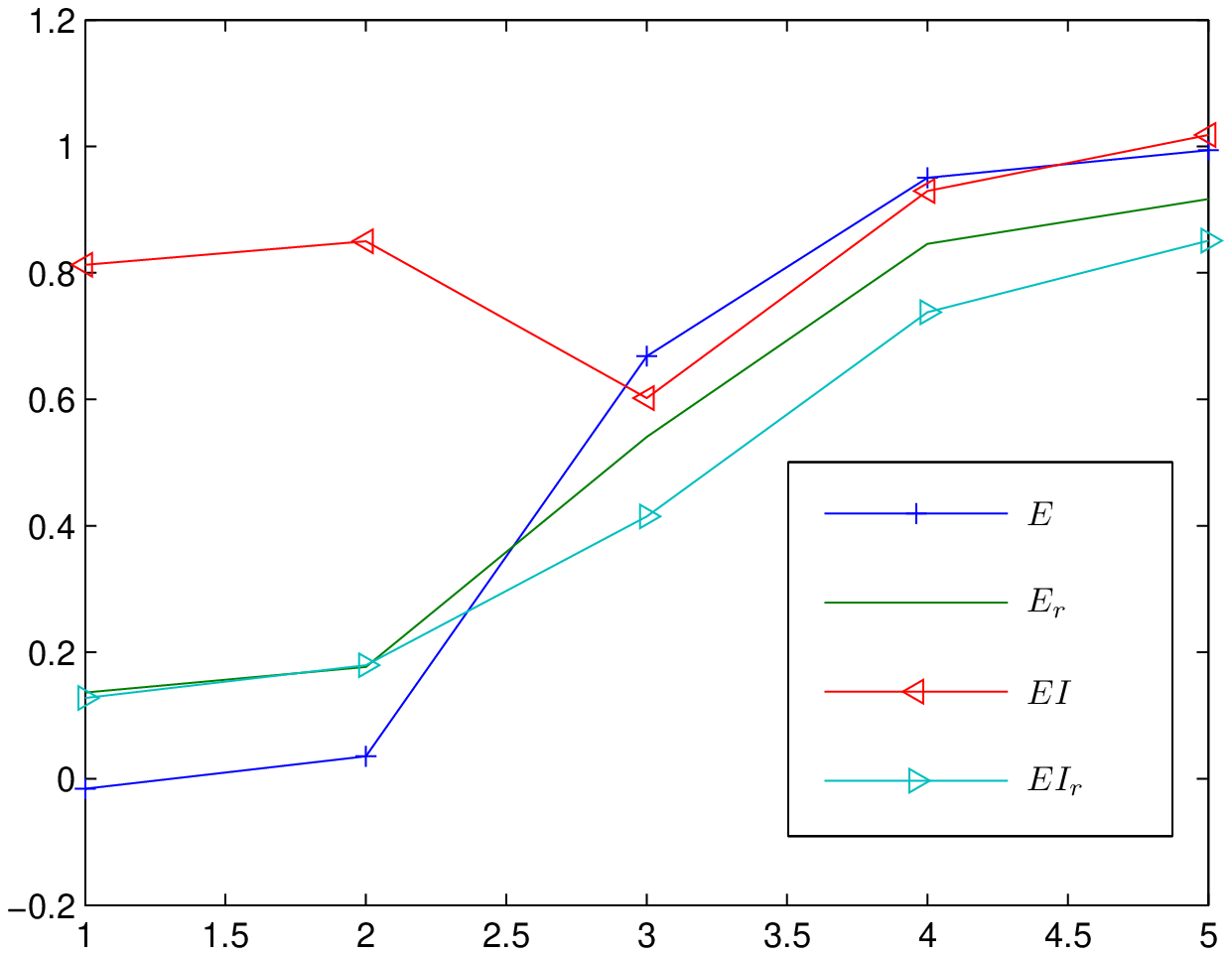}
  \put(-220,0){(a)}
  \put(0,15){\resizebox{5.8cm}{!}{\epsffile{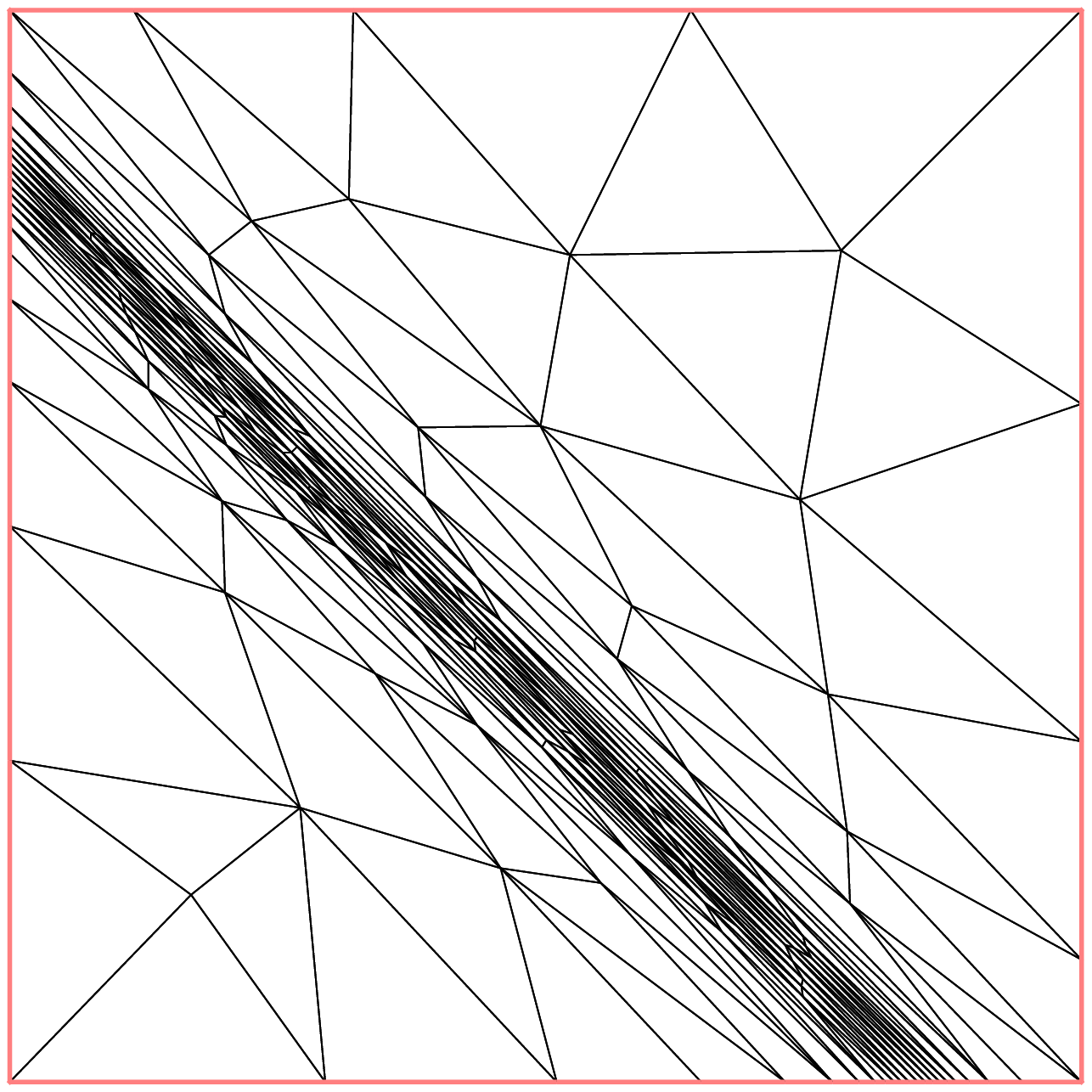}}}
  \put(-100,-200){\resizebox{9cm}{!}{\epsffile{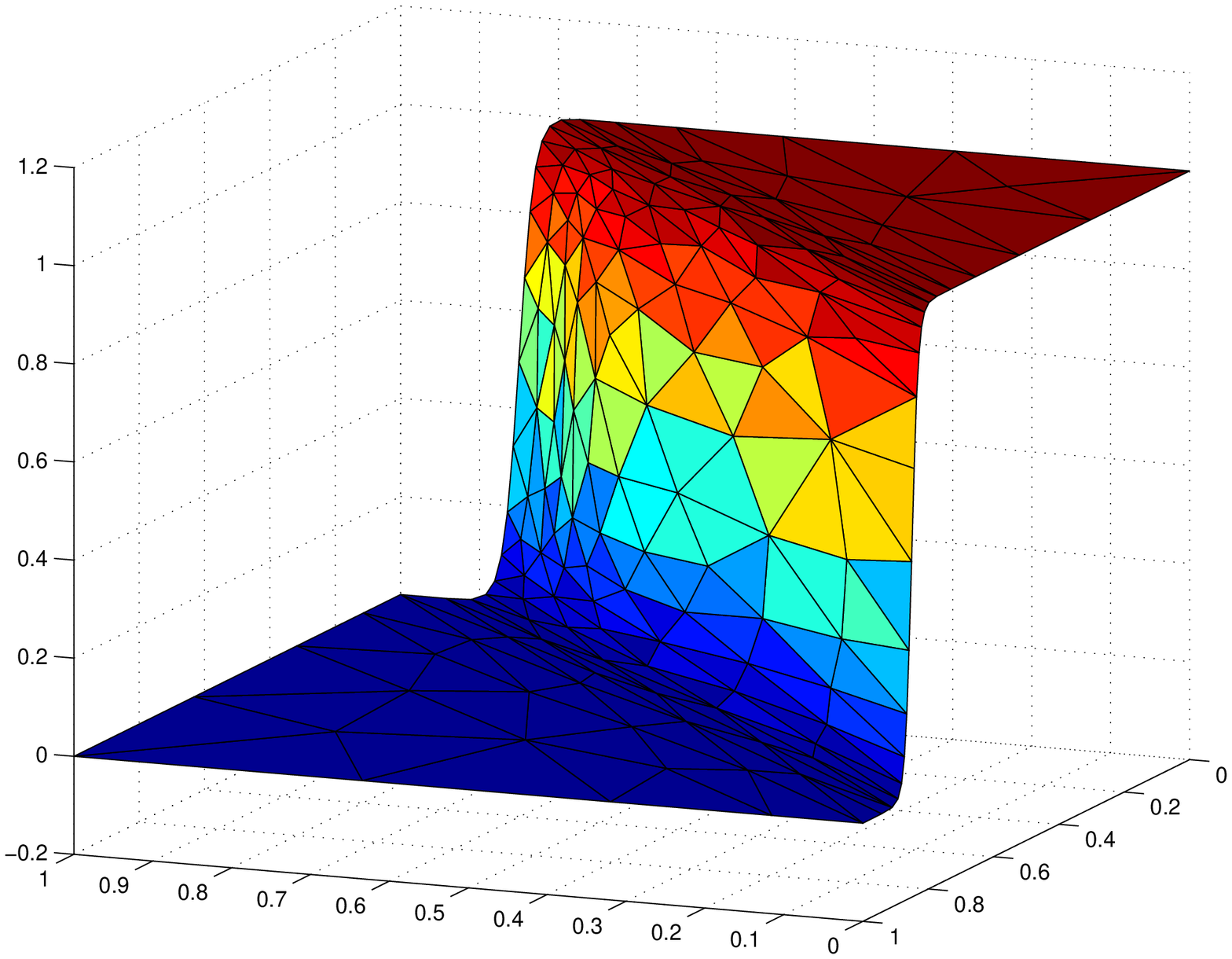}}}
 \put(10,0){(b)} \put(-100,-170){(c)}
 \put(-215,-200){\mbox{Figure 3: (a) Four estimators, (b) the final mesh and (c) $u_h$ of the example 4.2.}}
\end{figure}
\noindent {\bf Example 4.3} This example is to solve the boundary
value problem of Poisson's equation
\begin{eqnarray}
-\triangle u&=&f,\quad {\bf x}\in \Omega\equiv(0,1)\times(0,1),
\end{eqnarray}
where the Dirichlet boundary condition and the right-hand side term
are chosen such that the exact solution is given by
\begin{eqnarray}
u({\bf x})=e^{x_1^2-0.8}.
\end{eqnarray}
This is an extreme example for anisotropic behavior where the
function $u$ only depends on one variable $x_1$ or $x_2$. Such
functions are the real challenge in the a posteriori error analysis
since one is not allowed to use this knowledge. It is obvious our
four estimators perform very well. See Table 3 and Figure 4 for more
details.
\begin{center}
{\renewcommand{\baselinestretch}{1.0}\small
\begin{tabular}{cccccccc}
\multicolumn{8}{c}{\hbox{{\bf Table 3: Four estimators and $\delta$ in example~4.3}~~  }}\\
 \hline
step & $N$ &$E$ & $E_r$ &$EI$ &$EI_r$& $||H-H_r||$& $\delta$\\
\hline
1 & 26 & 0.863888 & 0.659135 & 1.14815  & 0.700543  & 1.02272 & -\\
2 & 26 & 0.899800 & 0.695822 & 1.13359  & 0.721021  & 0.976705 & -\\
3 & 32 & 0.984515 & 0.804626 & 1.13656  & 0.835906  & 0.808743 & 1.82\\
4 & 43 & 0.993768 & 0.938219 & 1.08240  & 0.967953  & 0.402943 & 4.72\\
5 & 66 & 0.990844 & 0.960663 & 1.04390  & 0.970201  & 0.190728 & 3.49\\
\hline
\end{tabular}}
\end{center}

\begin{figure}[ht]
  \includegraphics[width=8cm]{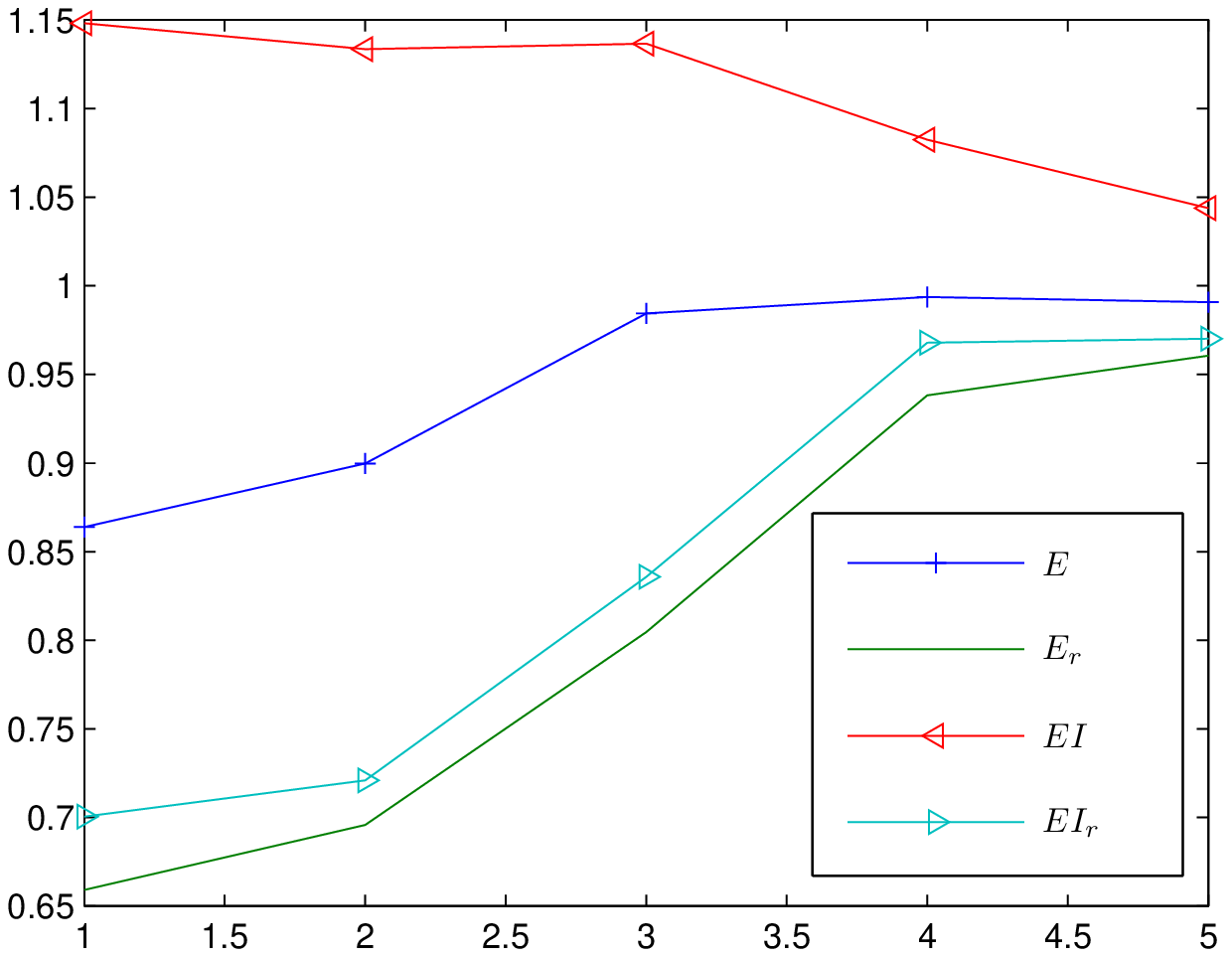}
  \put(-200,0){(a)}
  \put(-10,15){\resizebox{6.2cm}{!}{\epsffile{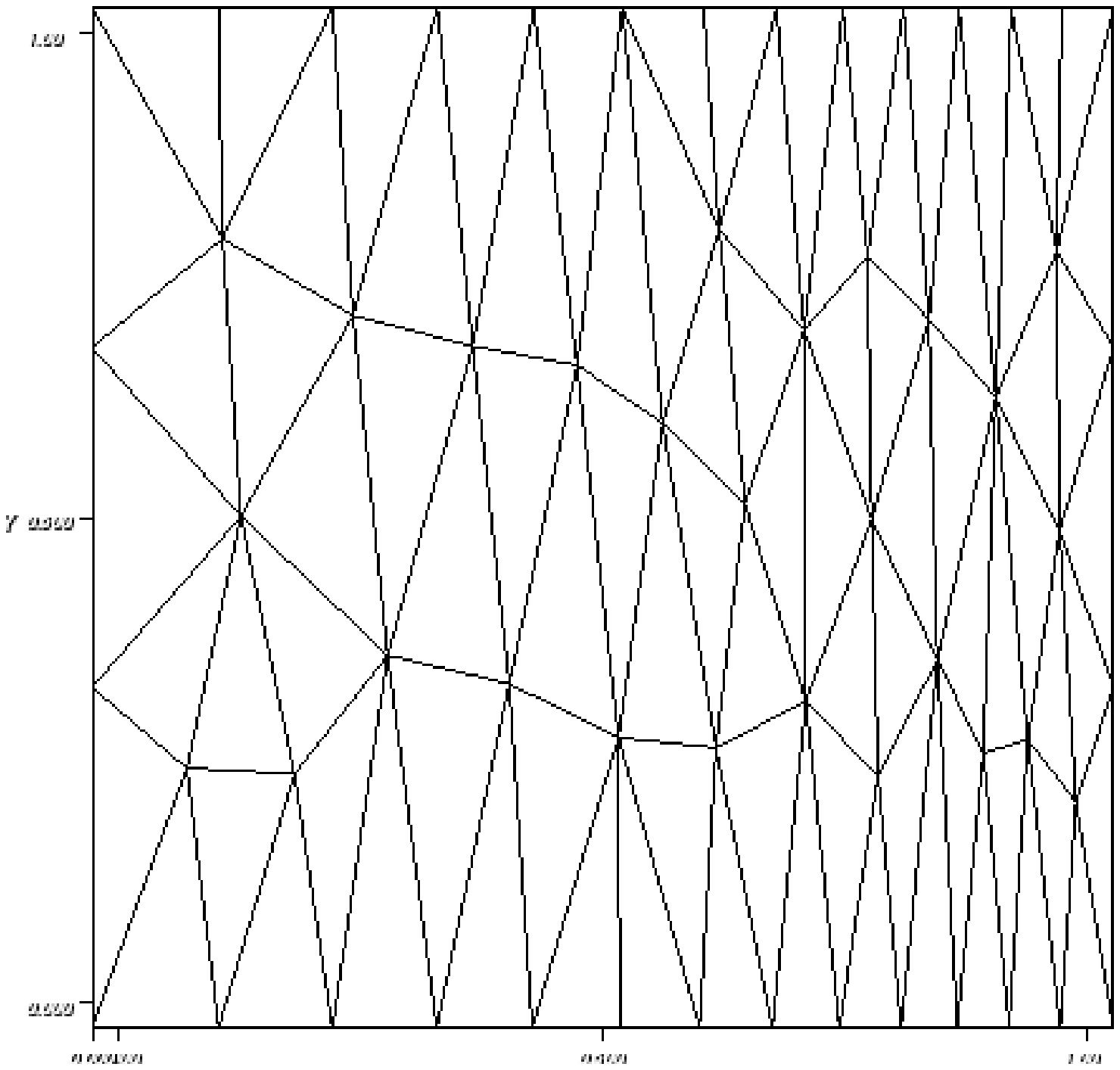}}}
  \put(10,0){(b)}\put(-90,-160){(c)}
    \put(-100,-200){\resizebox{9cm}{!}{\epsffile{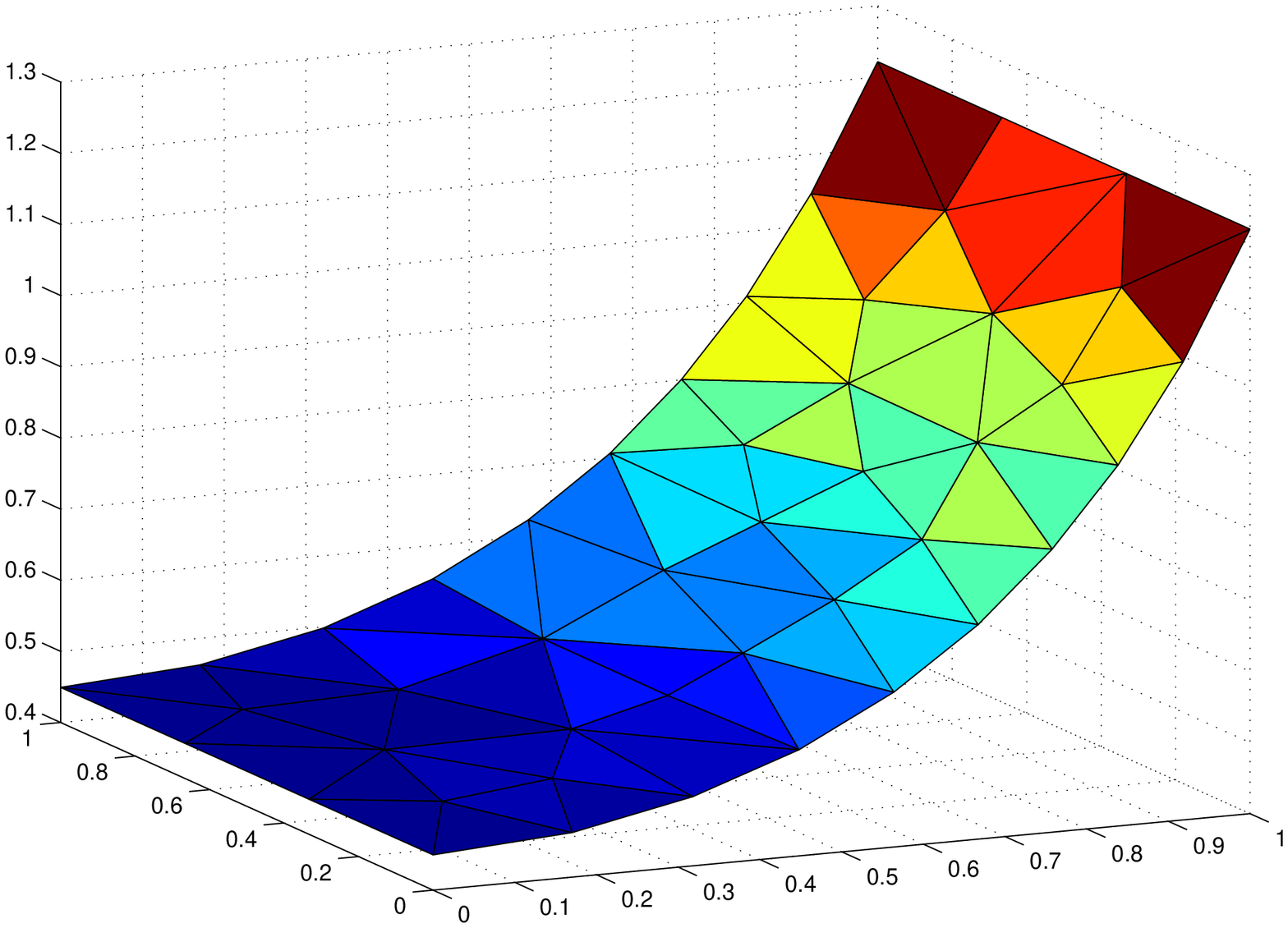}}}
 \put(-215,-200){\mbox{Figure 4: (a) Four estimators, (b) the final mesh and (c) $u_h$ of the example 4.3.}}
\end{figure}
 \noindent{\bf Example 4.4} This example is to solve the boundary value
problem of Poisson's equation
\begin{eqnarray}
-\triangle u&=&f,\quad {\bf x}\in \Omega\equiv(-1,1)\times(-1,1),
\end{eqnarray}
where the Dirichlet boundary condition and the right-hand side term
are chosen such that the exact solution is given by
\begin{eqnarray}
u({\bf x})=x_1^2x_2+x_2^3+\tanh(10(\sin(5x_2)-2x_1)).
\end{eqnarray}
The solution is anisotropic along the zigzag curve
$\sin(5x_2)-2x_1=0$ and changes sharply in the direction normal to
this curve(taken from \cite{Hecht2005,Lip}). For more details see
Table 4 and Figure 5.
\begin{center}
{\renewcommand{\baselinestretch}{1.0}\small
\begin{tabular}{cccccccc}
\multicolumn{8}{c}{\hbox{{\bf Table 4: Four estimators and $\delta$ in example~4.4}~~  }}\\
 \hline
step & $N$ &$E$ & $E_r$ &$EI$ &$EI_r$ & $||H-H_r||$& $\delta$\\
\hline
1 & 146 & -0.644460 & 0.106387 & 0.842915  & 0.082528  & 520.885 & -\\
2 & 325 & 0.014871 & 0.180764 & 0.534598  & 0.113211  & 490.853 & 0.15\\
3 & 756 & 0.642053 & 0.538530 & 0.559720  & 0.380420  & 289.914 & 1.25\\
4 & 1515 & 0.940350 & 0.854614 & 0.935262  & 0.766470  & 105.952 & 2.90\\
5 & 2826 & 0.993081 & 0.913258 & 1.08005  & 0.896111  &75.6472 & 1.08\\
\hline
\end{tabular}}
\end{center}

\begin{figure}[ht]
  \includegraphics[width=8cm]{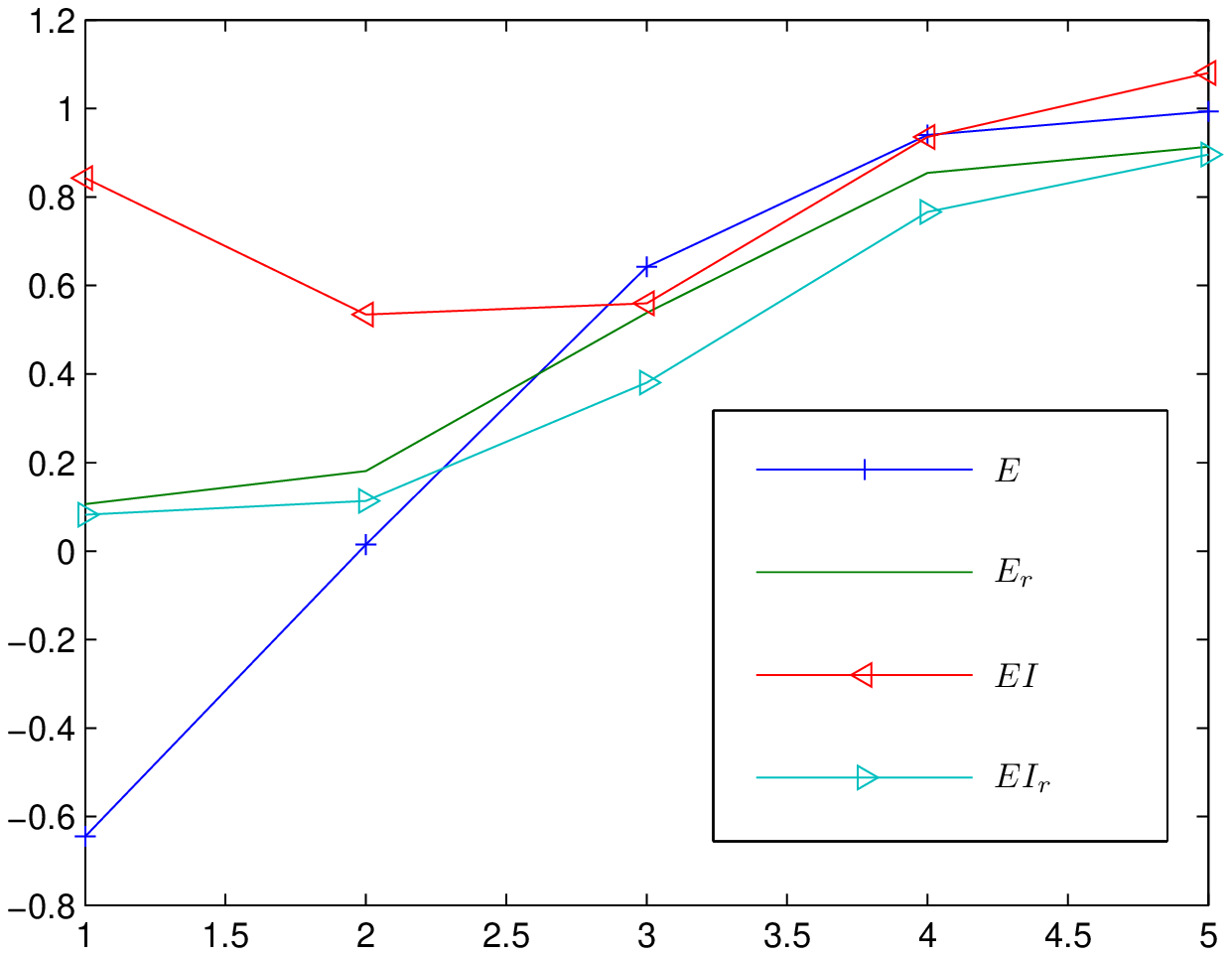}
  \put(-210,0){(a)}
  \put(-10,15){\resizebox{7.7cm}{!}{\epsffile{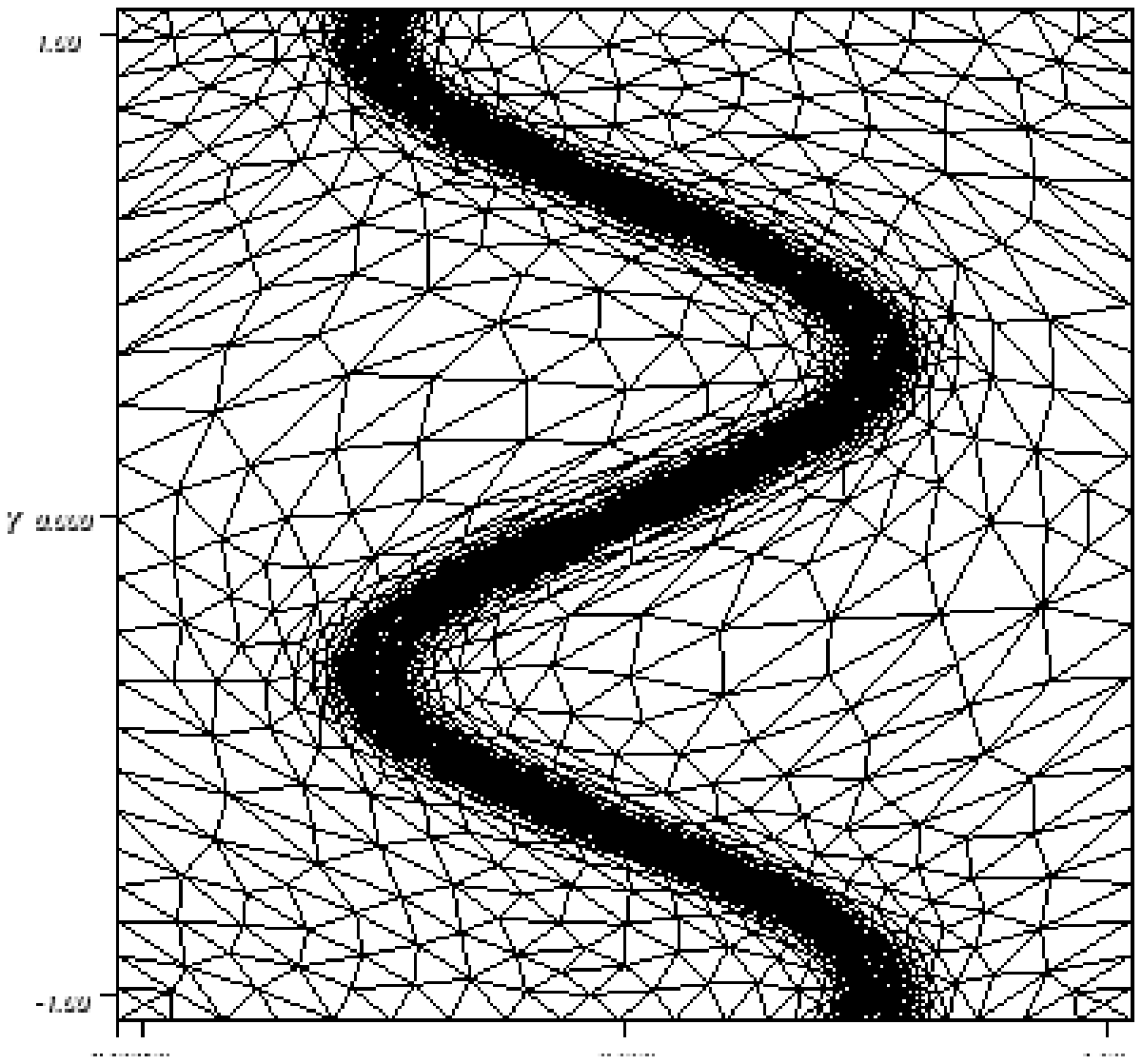}}}
\put(30,0){(b)} \put(-110,-170){(c)}
    \put(-100,-200){\resizebox{9cm}{!}{\epsffile{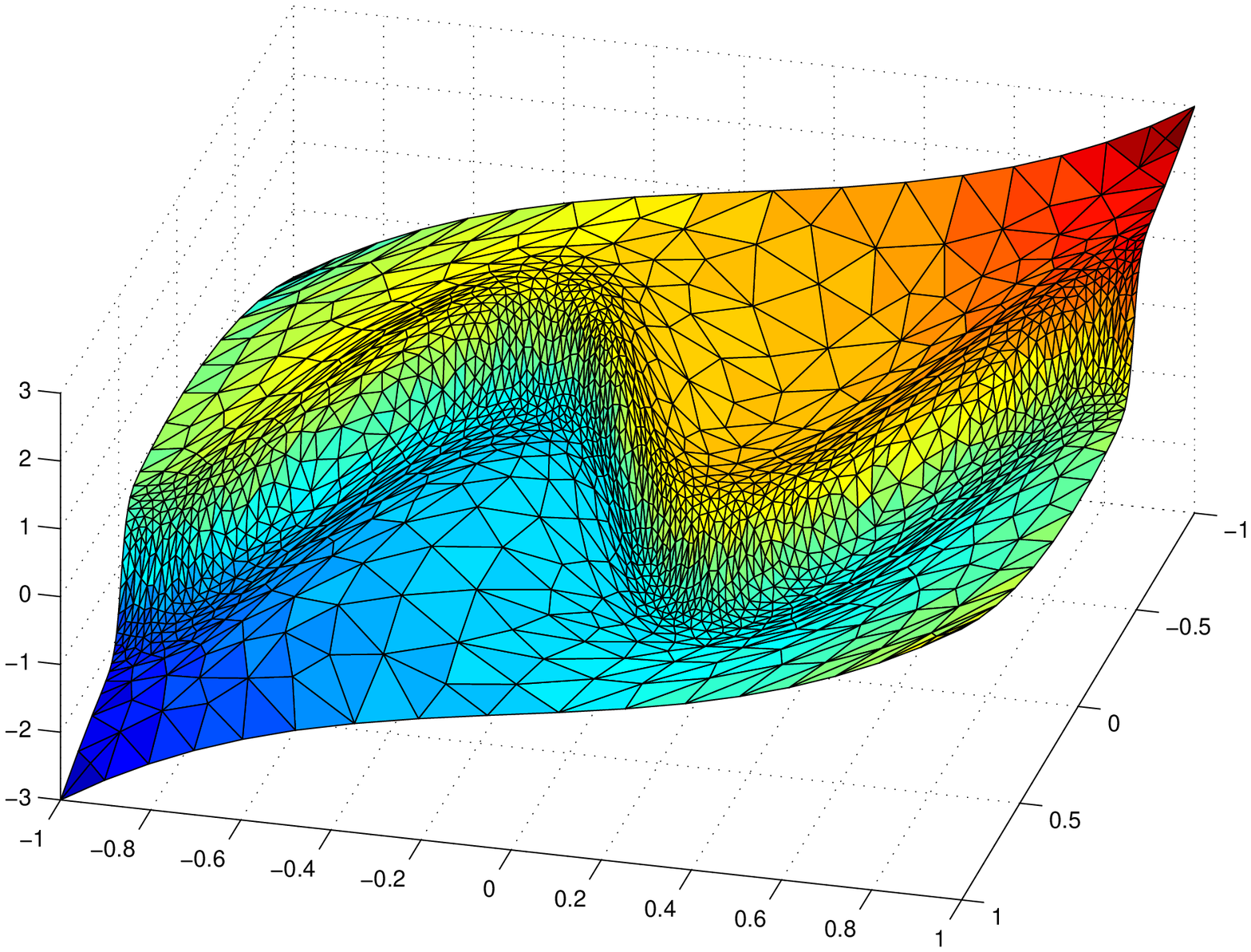}}}
 \put(-215,-210){\mbox{Figure 5: (a) Four estimators, (b) the final mesh and (c) $u_h$ of the example 4.4.}}
\end{figure}
\noindent{\bf Example 4.5} This example is to solve the boundary
value problem of
\begin{eqnarray}
-\triangle u&=&0,\quad {\bf x}\in
\Omega\equiv(-0.5,0.5)\times(0,0.5)\cup (-0.5,0)\times(-0.5,0).
\end{eqnarray}
The Dirichlet boundary condition is chosen such that the exact
solution is given by
$$u=r^{\frac{2}{3}}\sin(\frac{2}{3}\theta),$$
where $(r,\theta)\in \Omega$ are the usual polar coordinates. It is
well known that the exact solution $u\in
H^{\frac{5}{3}-\epsilon}(\Omega) (\forall \epsilon>0)$. So we expect
our estimators can be extended to more problems especially for those
with low regularity. See Figure 6 for more details.
\begin{figure}[ht]
  \includegraphics[width=8cm]{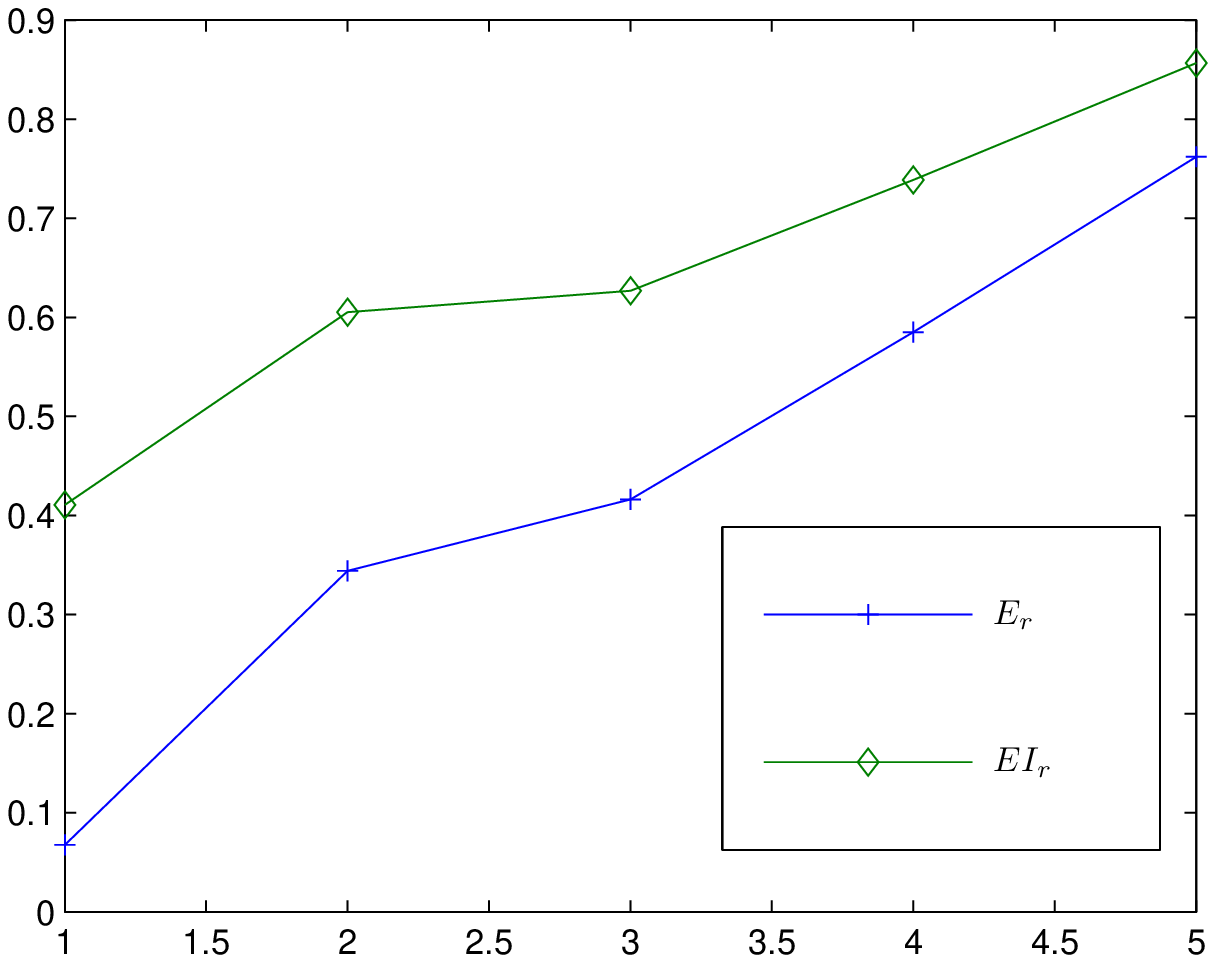}
  \put(-210,5){(a)}
  \put(-10,13){\resizebox{6.7cm}{!}{\epsffile{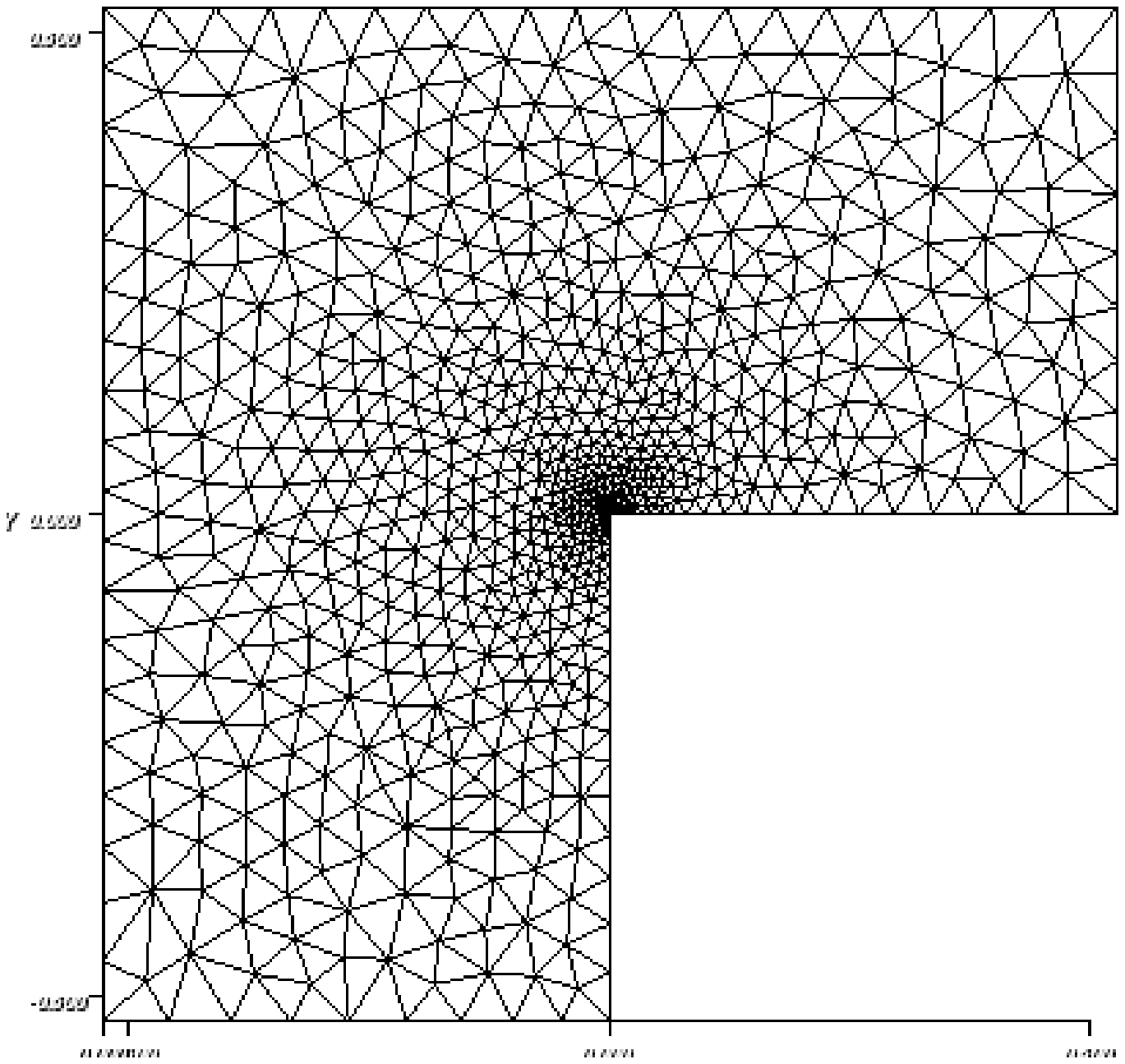}}}
\put(10,0){(b)} \put(-75,-180){(c)}
    \put(-100,-200){\resizebox{9cm}{!}{\epsffile{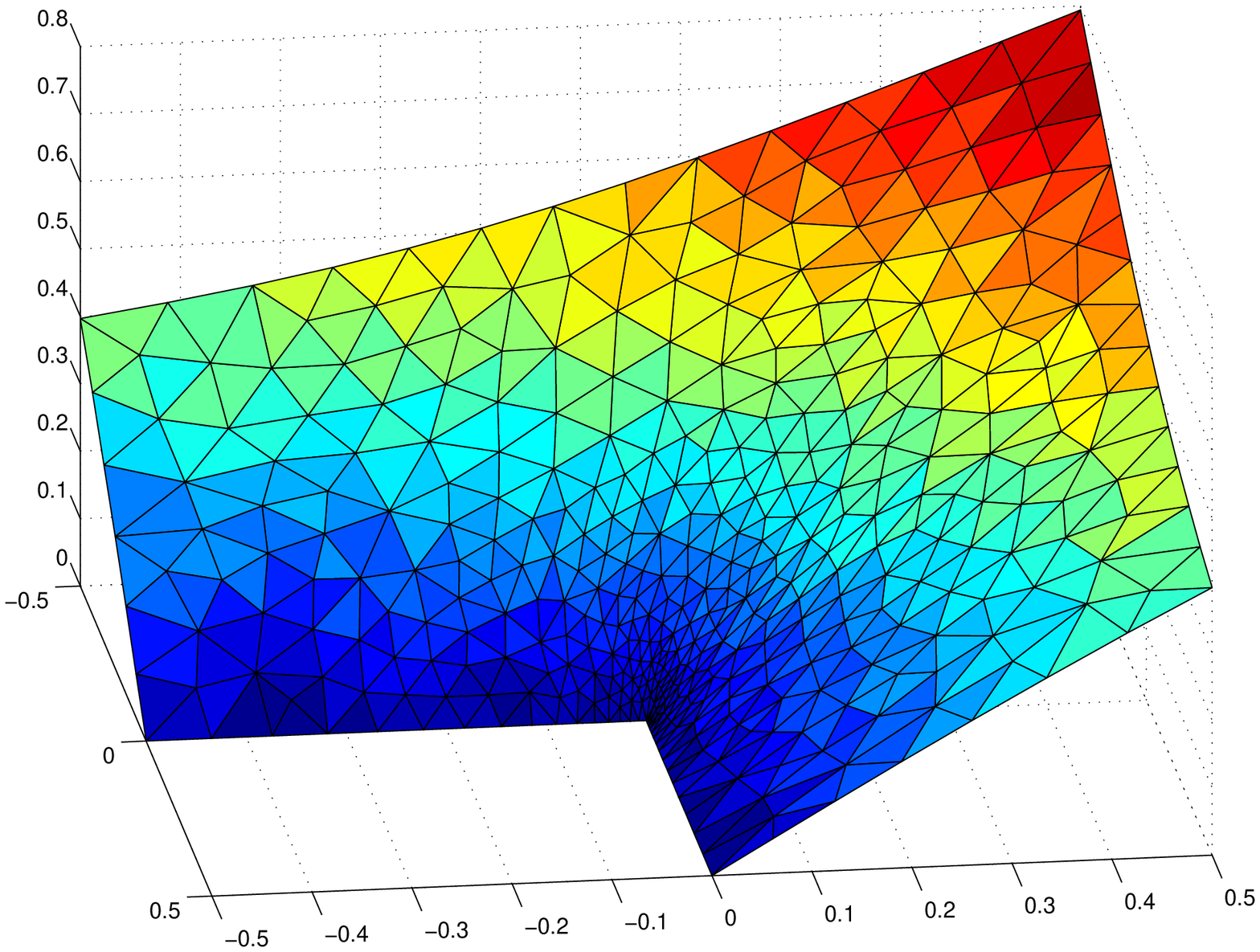}}}
 \put(-215,-210){\mbox{Figure 6: (a) Two estimators, (b) the final mesh and (c) $u_h$ of the example 4.5.}}
\end{figure}

\noindent{\bf Example 4.6} This example is to solve the boundary
value problem
\begin{eqnarray}
-\epsilon\triangle u+u=f,\quad {\bf x}\in
\Omega\equiv(0,1)\times(0,1),
\end{eqnarray}
where the Dirichlet boundary condition and the right-hand side term
are chosen such that the exact solution is the same as example
4.1(taken from \cite{Huang}). Note that
$$E=\frac{\eta^2}{\epsilon^{-1}\|u-u_h\|_{0,\Omega}^2+\|\nabla(u-u_h)\|_{0,\Omega}^2}$$
($E_r$,
$EI$ and $EI_r$ are defined similarly).
 See Table 6 and Figure
7 for more details.
\begin{figure}[ht]
  \centering
  \includegraphics[width=8.5cm]{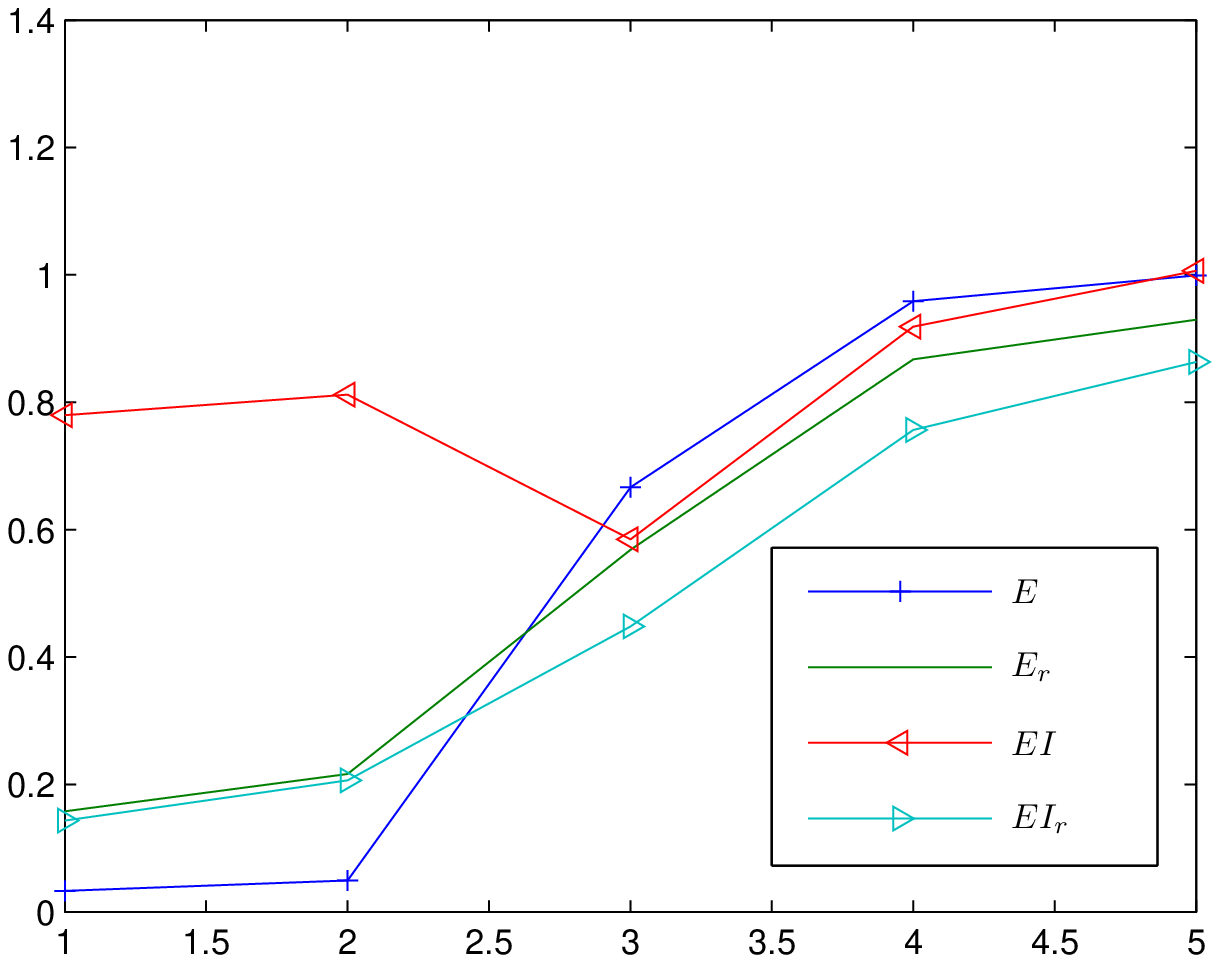}
  \put(-230,0){\mbox{Figure 7: Four estimators of example 4.6.}}
\end{figure}
\begin{center}
{\renewcommand{\baselinestretch}{1.0}\small
\begin{tabular}{cccccccc}
\multicolumn{8}{c}{\hbox{{\bf Table 6: Four estimators and $\delta$ in example~4.6}~~  }}\\
 \hline
step & $N$ &$E$ & $E_r$ &$EI$ &$EI_r$ & $||H-H_r||$& $\delta$\\
\hline
1 & 94  & 0.033119 & 0.158037 & 0.779738 & 0.143494 & 264.624 & -\\
2 & 113 & 0.049065 & 0.216479 & 0.811994 & 0.206575 & 229.060 & 1.57\\
3 & 189 & 0.666480 & 0.568354 & 0.584863 & 0.448080 & 148.367 & 1.69\\
4 & 272 & 0.958780 & 0.867167 & 0.918925 & 0.756519 & 58.9264 & 5.07\\
5 & 278 & 0.999348 & 0.929424 & 1.00642  & 0.863337 & 32.6764 & 54.05\\
\hline
\end{tabular}}
\end{center}
From experiments above we conclude that our a posteriori error
estimators $\eta$ and $\eta_I$ are always asymptotically exact under
various isotropic and anisotropic meshes. So we may guess that the
superapproximation always holds during the adaptive process.

\section{Conclusions}
In the previous sections we have developed a new type of a
posteriori error estimators suitable for moving mesh methods under
general meshes(especially anisotropic meshes). In our next paper we
want to design adaptive algorithms using the estimators, i.e., to
give a new metric tensor for moving mesh method.

\end{document}